\documentclass[conference,a4paper,11pt]{article}

\usepackage[utf8]{inputenc}
\usepackage[T1]{fontenc}
\usepackage{url}              
\usepackage{cite}             

\usepackage[cmex10]{amsmath}  
\usepackage{mleftright}       
\mleftright                   

\usepackage{graphicx}         
\usepackage{booktabs}         

\newtheorem{lemma}{Lemma}[section]

\newtheorem{theorem}[lemma]{Theorem}
\newtheorem{proposition}[lemma]{Proposition}
\newtheorem{corollary}[lemma]{Corollary}
\newtheorem{example}[lemma]{Example}

\newtheorem{remark}[lemma]{Remark}

\begin{document}

\title{Algebraic closures and their variations\footnote{The work
was carried out in the framework of the State Contract of the
Sobolev Institute of Mathematics, Project No.~FWNF-2022-0012, and
partially supported by Committee of Science in Education and
Science Ministry of the Republic of Kazakhstan, Grant No.
AP19677451.}}

  \author{Sergey V. Sudoplatov}

\date{}

\maketitle

\begin{abstract}
We study possibilities for algebraic closures, differences between
definable and algebraic closures in first-order structures, and
variations of these closures with respect to the bounds of
cardinalities of definable sets and given sets of formulae.
Characteristics for these possibilities and differences are
introduced and described. These characteristics are studied for
some natural classes of theories. Besides algebraic closure
operators with respect to sets of formulae are introduced and
studied. Semilattices and lattices for families of these operators
are introduced and characteristics of these structures are
described.
\end{abstract}

{\bf Key words:} algebraic closure, definable closure, degree of
algebraization, algebraic set, $\Delta$-algebraic element,
semilattice.

\section{Introduction}

The notions of definable and algebraic closures are broadly used
in Model Theory and applications \cite{Sh, Hodges, PiG, Marker,
TZ}.

In the paper, we study possibilities for algebraic closures,
differences between definable and algebraic closures in
first-order structures, and variations of these closures with
respect to the bounds of cardinalities of definable sets and given
sets of formulae. Characteristics for these possibilities and
differences are introduced and described. These characteristics
are studied for some natural classes of theories.

The paper is organized as follows. In Section 2, preliminary
notions, notations and properties for algebraic and definable
closures are considered. In Section 3, we study variations of
algebraic closure, properties and possibilities of these
variations, of degree of algebraization, and of difference between
definable and algebraic closures, illustrating these possibilities
by a series of examples. Algebraic sets and their degrees are
studied in Section 4. In Section 5, we consider algebraic closures
relative sets of formulae, study a hierarchy of operators of
algebraic closures relative these sets. We introduce semilattices
and lattices for families of these operators and describe some
characteristics of these structures.

Throughout we use both the standard model-theoretic terminology
and notions of Lattice Theory \cite{Birk, Grat, BS}.

\section{Pregeometries and closures}

{\bf Definition} \cite{Hodges, PiG, Marker, TZ}, cf. \cite{EmMe}.
For a set $S$, its Boolean $P(S)$, and an operator ${\rm
cl}\mbox{: }P(S)\to P(S)$, a pair $\langle S,{\rm cl} \rangle$ is
called a {\em pregeometry} or a {\em matroid}, if it satisfies the
following conditions for any $X\subseteq S$ and $a,b\in S$:

(i) $X\subseteq{\rm cl}(X)$ (Reflexivity);

(ii) ${\rm cl}({\rm cl}(X))={\rm cl}(X)$ (Transitivity);

(iii) if $a\in{\rm cl}(X)$ then $a\in{\rm cl}(Y)$ for some finite
$Y\subseteq X$ (Finite character).

(iv) if $a\in{\rm cl}(X\cup\{b\})\setminus{\rm cl}(X)$ then
$b\in{\rm cl}(X\cup\{a\})$ (Exchange property);

A pregeometry $\langle S,{\rm cl} \rangle$ is called a {\em
geometry} if:

(v) ${\rm cl}(\emptyset)=\emptyset$, and for every $a\in S$, ${\rm
cl}(\{a\})=\{a\}$ (ES-property).

The operator ${\rm cl}$ for the pregeometry is called the {\em
closure operator}. By the definition the closure operator for a
pregeometry produces a geometry if the closure of the empty set is
again empty and the closure of any singleton is again a singleton,
which, by reflexivity, equals that singleton.

\begin{remark} \label{pregeom}
{\rm \cite{PiG}. Any pregeometry $\langle S,{\rm cl} \rangle$
produces a {\em canonical} geometry $\langle S',{\rm cl}' \rangle$
putting $S'=\{{\rm cl}(\{a\})\mid a\in S\setminus{\rm
cl}(\emptyset)\}$ and for $X\subseteq S$, ${\rm cl}'(\{{\rm
cl}(\{a\})\mid a\in X\})=\{{\rm cl}(b)\mid b\in{\rm cl}(X)\}$.

In fact the ES-property for ${\rm cl}'$ is implies by the
reflexivity and transitivity for the operator~${\rm cl}$.}
\end{remark}

We fix a big saturated structure $\mathcal{M}$ and its theory
$T={\rm Th}(\mathcal{M})$. Following \cite{Spr} by a set $A$ of
$T$ we mean a subset of the universe $M$ in the structure
$\mathcal{M}$ satisfying some type ${\rm tp}(A,\mathcal{M})$.
Similarly, a tuple $\overline{a}$ in $T$ is a tuple of elements in
$\mathcal{M}$ satisfying given type ${\rm tp}(\overline{a})$.

In \cite{Sh}, two types of closures in structures are considered,
algebraic and definable, as well as the following concepts related
to them:

\medskip
{\bf Definition.} \cite{Sh, Hodges, TZ} 1. The tuple
$\overline{b}$ is {\em defined} by the formula
$\varphi(\overline{x},\overline{a})$ of $T$ with parameters
$\overline{a}$, if $\varphi(\overline{x},\overline{a})$ has unique
solution $\overline{b}$.

The tuple $\overline{b}$ is {\em defined} by the type $p$ if
$\overline{b}$ is the unique tuple which realizes $p$. It is {\em
definable} over a set $A$ if ${\rm tp}(\overline{b}/A)$ defines
it.

2. For a set $A$ of a theory $T$ the union of sets of solutions of
formulae $\varphi(x,\overline{a})$, $\overline{a}\in A$, such that
$\models\exists^{=n}x\:\varphi(x,\overline{a})$ for some
$n\in\omega$ (respectively
$\models$~$\exists^{=1}x\:\varphi(x,\overline{a})$) is said to be
an \emph{algebraic {\rm (}definable {\rm or} definitional{\rm )}
closure}\index{Closure!algebraic}\index{Closure!definable} of $A$.
An algebraic closure of $A$ is denoted by ${\rm
acl}(A)$\index{${\rm acl}(A)$} and its definable (definitional)
closure, by ${\rm dcl}(A)$\index{${\rm dcl}(A)$}.

In such a case we say that the formulae $\varphi(x,\overline{a})$
{\em witness} that algebraic / definable (definitional) closure,
and these formulae are called {\em algebraic / defining}.

Any element $b\in{\rm acl}(A)$ (respectively, $b\in{\rm dcl}(A)$)
is called {\em algebraic} ({\em definable} or {\em definitional})
over $A$. If the set $A$ is fixed or empty, we just say that $b$
is {\em algebraic}, ({\em definable}, or definitional).

3. If ${\rm dcl}(A)={\rm acl}(A)$, ${\rm cl}_1(A)$ denotes their
common value.

4. If $A={\rm acl}(A)$ (respectively, $A={\rm dcl}(A)$ ) then $A$
is called {\em algebraically} ({\em definably}) closed.

5. The type $p$ is {\em algebraic} ({\em defining}) if it is
realized by finitely many tuples (unique one) only, i.e., it
contains an algebraic (defining) formula $\varphi$. This formula
$\varphi$ can be chosen with the minimal number of solution, and
in such a case $\varphi$ isolates $p$. The number of these
solutions is called the {\em degree} ${\rm deg}(p)$ of $p$.

6. The complete algebraic types $p(x)\in S(A)$ are exactly ones of
the form ${\rm tp}(a/A)$, where $a$ is algebraic over $A$. The
{\em degree} of $a$ over $A$, ${\rm deg}(a/A)$ is the degree of
${\rm tp}(a/A)$.

\begin{remark} \label{acl-dcl-pregeom}
{\rm \cite{Sh}. The pairs $\langle M,{\rm acl}\rangle$ and
$\langle M,{\rm dcl}\rangle$ satisfy the following properties:

(i) the reflexivity: it is witnessed by the formula $x\approx y$;

(ii) the transitivity: if the formulae
$\varphi_1(x_1,\overline{a}),\ldots,\varphi_n(x_n,\overline{a})$
witnessed that $b_1,\ldots,b_n\in{\rm acl}(A)$ (respectively,
$b_1,\ldots,b_n\in{\rm dcl}(A)$) and the formula
$\psi(x,b_1,\ldots,b_n)$ witnesses that $c\in{\rm
acl}(\{b_1,\ldots,b_n\})$ (respectively, $c\in{\rm
dcl}(\{b_1,\ldots,b_n\})$) then the formula
\begin{equation}\label{eq_acl1}\exists
x_1,\ldots,x_n\left(\psi(x,x_1,\ldots,x_n)\wedge\bigwedge\limits_{i=1}^n\varphi_i(x_i,\overline{a})\right)\end{equation}
witnesses that $c\in{\rm acl}({\rm acl}(A))$ (respectively,
$c\in{\rm dcl}({\rm dcl}(A))$);

(iii) the finite character: if a formula $\varphi{x,\overline{a}}$
witnesses that $a\in{\rm acl}(A)$ (respectively, $a\in{\rm
dcl}(A)$) then $a\in{\rm acl}(A_0)$ for the finite $A_0\subseteq
A$ consisting of coordinates in $\overline{a}$.}
\end{remark}

\section{$n$-Closures and degrees of algebraization}

Now we consider modifications of the algebraic closure related to
the sets of solutions of formulas bounded in cardinality by some
natural number.

\medskip
{\bf Definition.} 1. For $n\in\omega\setminus\{0\}$ and a set $A$
an element $b$ is called {\em $n$-algebraic} over $A$, if
$a\in{\rm acl}(A)$ and it is witnessed by a formula
$\varphi(x,\overline{a})$, for $\overline{a}\in A$, with at most
$n$ solutions.

2. The set of all $n$-algebraic elements over $A$ is denoted by
${\rm acl}_n(A)$.

3. If $A={\rm acl}_n(A)$ then $A$ is called {\em
$n$-algebraically} closed.

4. The type $p$ is {\em $n$-algebraic} if it is realized by at
most $n$ tuples only, i.e., ${\rm deg}(p)\leq n$.

5. The complete $n$-algebraic types $p(x)\in S(A)$ are exactly
ones of the form ${\rm tp}(a/A)$, where $a$ is $n$-algebraic over
$A$, i.e., with ${\rm deg}(a/A)\leq n$. Here ${\rm deg}(a/A)=k\leq
n$ defines the {\em $n$-degree} ${\rm deg}_n(a/A)$ of ${\rm
tp}(a/A)$ and of $a$ over $A$.

6. If ${\rm acl}(A) = {\rm acl}_n(A)$ then minimal such $n$ is
called the {\em degree of algebraization} over the set $A$ and it
is denoted by ${\rm deg}_{\rm acl}(A)$. If that $n$ does not exist
then we put ${\rm deg}_{\rm acl}(A)=\infty$. The supremum of
values ${\rm deg}_{\rm acl}(A)$ with respect to all sets $A$ of
given theory $T$ is denoted by ${\rm deg}_{\rm acl}(T)$ and called
the {\em degree of algebraization} of the theory $T$.

7. Following \cite{Zilber} theories $T$ with ${\rm deg}_{\rm
acl}(T)=1$, i.e., with defined ${\rm cl}_1(A)$ for any set $A$ of
$T$, are called {\em quasi-Urbanik}, and the models $\mathcal{M}$
of $T$ are {\em quasi-Urbanik}, too.

\medskip
\begin{remark} \label{alg_type} {\rm By the definition any algebraic
type is $n$-algebraic for some $n$, it is isolated and have a
fixed finite set of solutions in any elementary extension of given
model. In particular, definable types are $1$-algebraic and have
unique realizations in any given model. We have ${\rm dcl}(A)={\rm
acl}_1(A)$ for any $A$.

Notice that any set $B$ of realizations of an algebraic type over
$A$ is the union of the subsets $B_n$ of realizations for its
$n$-algebraic subtypes for all $n$:
$A=\bigcup\limits_{n\in\omega}B_n$. Moreover, the sets $B_n$ form
an increasing chain, with $B_m=B_n$ for $m>n$.

Notice also that it is essential in the definition of
$n$-algebraic types that models for the sets of realizations are
arbitrary, since, for instance, a theory $T$ of infinitely many
disjoint nonempty unary predicates $U_i$, $i\in I$, has models
with arbitrarily finitely many realizations of nonisolated type
$p(x)=\{\neg U_i(x)\mid i\in I\}$. This type has infinitely many
realizations in an appropriate $T$-model, and therefore it is not
algebraic.}
\end{remark}

The following proposition allows to transform algebraic types into
their finite subtypes.

\begin{proposition}\label{prop_alg_1} {\rm (cf. \cite[Lemma 6.1]{Sh}).}
$1.$ A type $p$ defines a tuple if and only if for some finite
$q\subseteq p$, $\bigwedge q$ defines that tuple.

$2.$ A type $p$ is $n$-algebraic if and only if for some finite
$q\subseteq p$, $\bigwedge q$ is $n$-algebraic; moreover, $q$ can
be chosen with the same set of realizations as for $p$.
\end{proposition}

Recall \cite{Hodges} that for a set $A$ and an element $a$ the
{\em $A$-orbit} ${\rm Orb}_A(a)$ of $a$ is the set of all elements
$b$ in the given structure which are connected with $a$ by an
$A$-automorphism.

\medskip
The following proposition gives an algebraic characterization for
$n$-algebraic types.

\begin{proposition}\label{prop_alg_2} A type $p$ is $n$-algebraic over $A$ if and only if
any / some $(|A|+|T|)$-saturated model $\mathcal{M}$ containing
$A$ has finitely many $A$-orbits $O$ consisting of realizations of
$p$, all these orbits are finite, and moreover the union $\bigcup
O$ has at most $n$-elements. If $p$ is complete then that
$A$-orbit is unique in $\mathcal{M}$.
\end{proposition}

Proof. Taking a $(|A|+|T|)$-saturated model $\mathcal{M}$ with
$A\subseteq M$ and a $n$-algebraic type $p$ over $A$ we have
finitely many possibilities for links, by $A$-automorphisms,
between realizations of $p$. Moreover, if $b_1,b_2,\ldots,b_n$ are
all realizations of $p$ then for any $b_i$, $b_j$ either ${\rm
tp}(b_i/A)\ne{\rm tp}(b_j/A)$ and $b_i$, $b_j$ belong to distinct
$A$-orbits, or ${\rm tp}(b_i/A)={\rm tp}(b_j/A)$ and $b_i$, $b_j$
belong to a common $A$-orbit connected by an $A$-automorphism $f$
with $f(b_i)=f(b_j)$. Additionally, if $c$ is not a realization of
$p$ then $c$ does not belong to orbits of elements
$b_1,b_2,\ldots,b_n$.

Conversely, if any / some $(|A|+|T|)$-saturated model
$\mathcal{M}$ containing $A$ has finitely many $A$-orbits $O$
consisting of realizations of $p$ and the union $\bigcup O$ has at
most $n$-elements then by the saturation of $\mathcal{M}$ the type
$p$ can not have more than $n$ realizations implying that $p$ is
$n$-algebraic over $A$.

\begin{remark} \label{alg_type2} {\rm
Taking cycles of arbitrarily large diameters we can find an
algebraic type with arbitrary large and arbitrary many finite
orbits. At the same time, it is essential in Proposition
\ref{prop_alg_2} that the model is saturated, since following
Remark \ref{alg_type} there are non-algebraic types with
arbitrarily many realizations, and taking copies of Example in
Remark \ref{alg_type} we can obtain arbitrarily many finite
orbits.

Notice also that the properties of $n$-algebraicity over a set $A$
are both syntactic and semantic, since on the one hand it is
written by a family of formulae
$\varphi(x,\overline{a})\wedge\exists^{\leq
n}x\varphi(x,\overline{a})$, for $\overline{a}\in A$, and on the
other hand it is checked in any model $\mathcal{M}$ of given
theory $T$ with $A\subseteq M$ and satisfying the given type ${\rm
tp}(A)$.

For instance, taking a structure $\mathcal{M}$ with $A\subseteq M$
such that each element of $M$ is the unique solution of an
appropriate formula $\varphi(x,\overline{a})$, for
$\overline{a}\in A$, then $M$ is definably closed, i.e.,
$1$-algebraically closed, with unique automorphism $f$ fixing $A$
pointwise and this automorphism is identical. Following
\cite{Hodges} the structure $\mathcal{M}$ with unique automorphism
is called {\em rigid}. Thus the $1$-algebraicity of each element
in $M$ over the empty set, i.e., the condition $M={\rm
dcl}(\emptyset)$ produces the rigid structure $\mathcal{M}$.

We separate the forms of rigidity of a structure as follows: the
{\em semantic} one is defined in terms of trivial automorphism
group, and the {\em syntactic} one is in terms of $1$-algebraicity
of the universe over the empty set.

Thus, any syntactically rigid structure is semantically rigid and
it is defined uniquely. But not vice versa, in general. Indeed,
taking a structure $\mathcal{M}$ consisting of infinitely many
distinct constants $c_i$, $i\in I$, we have both the semantic and
syntactic rigidity. Extending $\mathcal{M}$ by a single
realization of the nonisolated type $p(x)=\{\neg x\approx c_i\mid
i\in I\}$ we obtain a structure $\mathcal{N}\succ\mathcal{M}$
preserving the semantic rigidity but loosing the syntactic one.
Taking disjoint unions \cite{Wo} of copies of $\mathcal{M}$ one
can obtain semantically rigid elementary extensions with
arbitrarily many new elements. Hence, in appropriate cases, the
semantic rigidity can be transformed for elementary extensions.

These appropriate cases admits that one-element finite orbits are
preserved under elementary extensions. We observe this
preservation for structures whose universes consist of constants
only. And one-element finite orbits are not preserved if, for
instance, a structure consists of infinitely many two-element
equivalence classes expanded by constants for all elements of
these classes.

Similarly, following Proposition \ref{prop_alg_2} homogeneous
structures $\mathcal{M}$ with at most $n$-element finite orbits,
under appropriate additional conditions, produce the
$n$-algebraicity. And some opposite conditions imply the negation
of $n$-algebraicity for $\mathcal{N}\equiv\mathcal{M}$.

Thus, in general case, poor automorphism groups, in particular,
identical automorphism groups ${\rm Aut}(\mathcal{M})$ can not
reflect adequately links for algebraic closures in structures
$\mathcal{N}\equiv\mathcal{M}$.}
\end{remark}

Recall that a theory $T$ is {\em $n$-transitive}, for
$n\in\omega$, if the type ${\rm tp}(\overline{a})$ of each tuple
$\overline{a}$ with $l(\overline{a})=n$ is forced by formulae
describing that coordinates of $\overline{a}$ coincide or do not
coincide.

\begin{remark} \label{alg_type3} {\rm
Clearly, any $n$-transitive theory is $m$-transitive for every
$m\leq n$. Besides, algebraic characteristics for $n$-transitive
theories such as degrees of algebraicity are really defined by
sets $A$ of cardinalities $\geq n$, since ${\rm acl}(B)=B$ for all
smaller sets $B$.

Thus, in general case, structures for ${\rm acl}(B)$ with small
finite cardinalities $|B|$ does not reflect general behavior of
the algebraic closure operator. It is valid for the operators
${\rm acl}_n(\cdot)$, too.}
\end{remark}

\begin{remark} \label{alg_type4} {\rm
In view of Remarks \ref{alg_type2} and \ref{alg_type3}, nir poor
automorphism groups, for structure similar to rigid, nor rich
automorphism groups, for structure similar to transitive, do not
reflect the frame of given structures and the frame for their
theories in the level of algebraic closure. }
\end{remark}

\begin{proposition}\label{prop_alg_3} {\rm (cf. \cite[Lemma 6.2]{Sh}).}
$1.$ $A\subseteq{\rm acl}_m(A)\subseteq{\rm acl}_n(A)\subseteq{\rm
acl}(A)$, for any $0<m\leq n$. In particular, $\langle M,{\rm
acl}_n\rangle$ is reflexive for any $n$.

$2.$ If $A\subseteq B$ and $n\geq 1$ then ${\rm
acl}_n(A)\subseteq{\rm acl}_n(B)$.

$3.$ If $A$ is definably {\rm (}algebraically{\rm )} closed then
$A={\rm dcl}(A)$ {\rm (}$A={\rm acl}(A)${\rm )}.

$4.$ If $A$ is $n$-algebraically closed for $n\geq 2$ then $A={\rm
acl}(A)$ iff any finite orbit over $A$ has at most $n$ elements.

$5.$ A tuple $\overline{b}$ is definable {\rm (}algebraic{\rm )}
over $A$ if and only if $\overline{b}\in{\rm dcl}(A)$ {\rm
(}$\overline{b}\in{\rm acl}(A)${\rm )}.
\end{proposition}

Proof. Items 1, 2, 3, 5 immediately follow by the definition and
\cite[Lemma 6.2]{Sh}. Item 4 is implied by Proposition
\ref{prop_alg_2}.

\medskip
Proposition \ref{prop_alg_3} immediately implies:

\begin{corollary}
If ${\rm cl}_1(A)$ exists then it equals ${\rm acl}_n(A)$ and it
is $n$-algebraically closed for each~$n$.
\end{corollary}

\begin{corollary}
For any structure $\mathcal{M}$, the pairs $\mathcal{S}_n=\langle
M,{\rm acl}_n\rangle$, $n\in\omega$, define ascending chains
$({\rm acl}_n(A))_{n\in\omega}$ for each $A\subseteq M$, where
${\rm acl}_0(A)=A$. These pairs $S_n$ coincide starting with some
$n_0$ iff all chains $({\rm acl}_n(A))_{n\geq n_0}$ are
singletons, and iff any finite orbit over arbitrary $A\subseteq M$
has at most $n_0$ elements in a homogeneous elementary extension
of $\mathcal{M}$.
\end{corollary}

\begin{remark} \label{rem_bound}
{\rm Since for $k\leq m$ there are $A^k_m=\frac{m!}{(m-k)!}$
$k$-tuples consisting of pairwise distinct elements in the given
$m$-element set we can not assert Item 5 of Proposition
\ref{prop_alg_3} for $n$-algebraic tuples. Indeed, orbits for
tuples can be much more greater than orbits for elements. In
particular, it holds for structures in the empty language.}
\end{remark}

\begin{remark} \label{rem_refl_fin_char}
{\rm Notice that along with the reflexivity for ${\rm acl}_n$, by
the definition any operator ${\rm acl}_n$ has the finite
character. At the same time ${\rm acl}_n$ can be non-transitive.
Indeed, taking an element $a$ and ${\rm acl}_n({\rm
acl}_n(\{a\}))$ one can find a structure $\mathcal{M}$, say a
graph, with $n$ new elements $b_1,\ldots,b_n$ as solutions of
$n$-algebraic formula $\varphi(x,a)$ isolating an $n$-algebraic
type $p(x)$ such that the $n$-algebraic formulae
$\varphi(x,b_1),\ldots,\varphi(x,b_n)$ produce $n^2$ new elements
$c_1,\ldots,c_{n^2}$ connected by $\{a\}$-orbits. Thus,
$c_1,\ldots,c_{n^2}\in {\rm acl}_n({\rm acl}_n(\{a\}))$ whereas
these elements do not belong to ${\rm acl}_n(\{a\})$.}
\end{remark}

\begin{remark} \label{rem_pregeom_n}
{\rm Clearly, Exchange property for ${\rm acl}_n(\cdot)$ can fail
even if ${\rm acl}(\cdot)$ satisfies it. Indeed, let $\Gamma$ be a
bipartite graph with parts $U$ and $V$ such that $|U|=n$,
$|V|=m>n$ and each vertex in $U$ is connected by an edge with each
vertex in $V$. We have $a\in{\rm acl}_n(\{b\})$ for any $a\in U$
and $b\in V$, whereas $b\in{\rm acl}(\{a\})\setminus{\rm
acl}_n(\{a\})$.

At the same time replacing ${\rm acl}_n(\cdot)$ by ${\rm
acl}_m(\cdot)$ we obtain that Exchange property. In fact, Exchange
property for ${\rm acl}_m(\cdot)$ is implied by Exchange property
for the case ${\rm acl}(\cdot)={\rm acl}_m(\cdot)$.}
\end{remark}

\begin{remark} \label{rem_equiv}
{\rm Let $T$ be a theory of an equivalence relation $E$. By the
definition any element $a$ of a finite $E$-class $X$ is
$|X|$-algebraic over any singleton $\{b\}\in X\setminus\{a\}$.
Moreover, $a$ is $n$-algebraic over $\emptyset$ if there are
finitely many $E$-classes of the cardinality $|X|$. Otherwise, if
there are infinitely many $E$-classes of the cardinality $|X|$
these $E$-classes do not contain algebraic elements over
$\emptyset$.

If the $E$-class $X$ is infinite it does not contain algebraic
elements and for each $A\subseteq X$, then ${\rm acl}(A)={\rm
dcl}(A)=A$, and the correspondent theory is quasi-Urbanik.}
\end{remark}

By Remark \ref{rem_equiv} we immediately obtain the following
dichotomy for values ${\rm deg}_{\rm acl}(T)$ of theories $T$ in
the language $\{E\}$ of an equivalence relation:

\begin{proposition}\label{deg_acl_equiv}
For any theory $T$ of an equivalence relation $E$ and a set $A$
either ${\rm deg}_{\rm acl}(A)\in\omega$, if $A$ consists of
elements with bounded orbits inside finite $E$-classes and between
these $E$-classes, or ${\rm deg}_{\rm acl}(A)=\infty$, otherwise.
\end{proposition}

\begin{example}\label{ex1}
{\rm Let $\mathcal{M}$ be a structure in the language $\{E\}$ of
equivalence relation and consists of $E$-classes with unbounded
distinct finite cardinalities then the sets ${\rm
acl}_n(\emptyset)$, $n\geq 1$, form an unboundedly increasing
sequence of $n$-algebraically closed sets. We observe the same
effect for ${\rm acl}_n(\emptyset)$, $n\geq 1$, where $A$ consists
of boundedly many elements in each $E$-class producing ${\rm
deg}_{\rm acl}(A)=\infty$ and ${\rm deg}_{\rm acl}({\rm
Th}(\mathcal{M}))=\infty$ as well.}
\end{example}

The following example illustrates the possibility ${\rm deg}_{\rm
acl}(T)=\omega$.

\begin{example}\label{ex2}
{\rm Let $\mathcal{M}$ be a structure in the language $\{E\}$ and
consisting of an equivalence relation $E$ with infinitely many
$n$-element $E$-classes for each $n\in\omega\setminus\{0\}$. By
the definition all finite orbits over a set $A\subseteq M$ are
exhausted by automorphisms transforming elements in $E$-classes
$X\subseteq M$ with $X\cap A\ne\emptyset$. By the arguments in
Example \ref{ex1} we have ${\rm deg}_{\rm acl}(T)=\infty$ for
$T={\rm Th}(\mathcal{M})$.

Now we expand the structure $\mathcal{M}$ till a structure
$\mathcal{M}'$ by countably many ternary predicates $R$ in the
following way:

1) if $\models R(a,b,c)$ then $a$ and $b$ belong to disjoint
$E$-classes $E(a)$ and $E(b)$ and $c\in E(a)\cup E(b)$ is the
unique solution of $R(a,b,z)$, moreover, we require that each
$c'\in E(a)\cup E(b)$ has a ternary symbol $R'$ with $\models
R'(a,b,c')\wedge\exists^{=1}z R'(a,b,z)$;

2) $\mathcal{M}'$ has same orbits as $\mathcal{M}$ over subsets of
$E$-classes and one-element orbits only in $E(a)\cup E(b)$, for
$E(a)\cap E(b)=\emptyset$, over sets $A\subseteq E(a)\cup E(b)$
which have both elements in $E(a)$ and in $E(b)$.

The structure $\mathcal{M}'$ can be formed step-by-step using an
appropriate generic construction \cite{CCMCT, Su073}.

By the definition for each $A\subseteq M'$, ${\rm deg}_{\rm
acl}(A)$ is finite, moreover, ${\rm deg}_{\rm acl}(A)=1$ for $A$
laying in several finite $E$-classes. At the same time the values
${\rm deg}_{\rm acl}(A)$ are unbounded in $\omega$ using subsets
of finite $E$-classes. Therefore, ${\rm deg}_{\rm acl}({\rm
Th}(\mathcal{M}'))=\omega$.}
\end{example}

\begin{theorem}\label{deg_acl}
$1.$ For any consistent theory $T$, ${\rm deg}_{\rm
acl}(T)\in((\omega+1)\setminus\{0\})\cup\{\infty\}$.

$2.$ For any $\lambda\in((\omega+1)\setminus\{0\})\cup\{\infty\}$
there is a theory $T_\lambda$ such that  ${\rm deg}_{\rm
acl}(T_\lambda)=\lambda$.
\end{theorem}

Proof. Item 1 holds since ${\rm deg}_{\rm acl}(T)$ is a supremum
of values in $(\omega\setminus\{0\})\cup\{\infty\}$.

2. At first we consider theories $T$ of equivalence relations $E$.
The values ${\rm deg}_{\rm acl}(T_n)=n\in(\omega\setminus\{0\})$
are realized by $E$ on $n$-element sets $M$ with unique
equivalence classes. Here each set $A\subseteq M$ has ${\rm
acl}(A)={\rm acl}_n(A)=M$. We obtain a similar effect adding
infinite $E$-classes, or taking finitely many finite $E$-classes
having a total of $n$ elements, or adding both these $E$-classes.

The value ${\rm deg}_{\rm acl}(T_\infty)=\infty$ is confirmed by
an equivalence relation $E$ with infinitely many finite
$E$-classes with unbounded finite cardinalities. Here we take a
set $A$ containing elements in all these classes producing ${\rm
acl}(A)\ne{\rm acl}_n(A)$ for any $n\in\omega$.

The value ${\rm deg}_{\rm acl}(T_\omega)=\omega$ is witnessed by
Example \ref{ex2}.

\begin{example}\label{ex3} {\rm By the definition any definable element is $n$-algebraic,
for any $n\geq 1$, but not vice versa. For instance, taking a
graph $\Gamma=\langle\{a_1,a_2\};\{e\}\rangle$, where $e$ is the
edge $[a_1,a_2]$, we obtain ${\rm acl}(\emptyset)=\{a_1,a_2\}$
witnessed by the formula $x\approx x$, i.e., ${\rm deg}_{\rm
acl}(\emptyset)={\rm deg}_{\rm acl}({\rm Th}(\Gamma))=2$, whereas
${\rm dcl}(\emptyset)=\emptyset$, since $a_1$ and $a_2$ are
connected by an automorphism. It is a minimal example, by
inclusion, in a graph language.

We observe the same effect for a structure with at least two
elements in the empty language, or in the graph language with the
empty relation, or in the language of one unary function
consisting of loops or of a cycle.}
\end{example}

\begin{example}\label{ex4} {\rm
For any linearly ordered structure $\mathcal{M}$ and a set
$A\subseteq M$, \begin{equation}\label{eq_lin}{\rm acl}(A)={\rm
dcl}(A),\end{equation} since for any formula
$\varphi(x,\overline{a})$ with finitely many solutions the finite
set $B$ of these solutions is linearly ordered and each element of
$B$ is defined by its position in that finite linear order. Thus,
for the linearly ordered structure $\mathcal{M}$, ${\rm deg}_{\rm
acl}(\mathcal{M})=1$, and it has a quasi-Urbanik theory.

Example 2.12 in \cite{KulMac} illustrates that for the circularly
ordered structure $c(\omega+\omega^\ast+{\bf
Q}+\omega+\omega^\ast+{\bf Q})$, ${\rm acl}(\emptyset)$ is a
proper superset of ${\rm dcl}(\emptyset)$ contrasting with the
linearly ordered case, since the first elements of the copies of
$\omega$ are connected by an automorphism, and the last elements
of the copies of $\omega^\ast$ are connected by an automorphism,
too.

At the same time, for any circularly ordered structure
$\mathcal{M}$ and a nonempty set $A\subseteq M$, ${\rm
acl}(A)={\rm dcl}(A)$, since any element $a\in A$ defines a linear
order on $M\setminus\{a\}$ satisfying the equality (\ref{eq_lin}).
Thus, for the circularly ordered structure $\mathcal{M}$, ${\rm
deg}_{\rm acl}({\rm Th}(\mathcal{M}))={\rm deg}_{\rm
acl}(\emptyset)$. Adding copies of $\omega+\omega^\ast+{\bf Q}$ to
the circularly ordered structure $c(\omega+\omega^\ast+{\bf
Q}+\omega+\omega^\ast+{\bf Q})$ we can obtain an arbitrary natural
value ${\rm deg}_{\rm acl}(\emptyset)$ for a circularly ordered
theory.

Continuing the process with infinitely and densely many copies of
$\omega+\omega^\ast+{\bf Q}$ and marking some of them circularly
by unary predicates we can obtain a circularly ordered theory $T$
with \begin{equation}\label{deg_inf}{\rm deg}_{\rm acl}(T)={\rm
deg}_{\rm acl}(\emptyset)=\infty.\end{equation} We can, for
instance, choose a prime number $p$ and mark $p$ copies
$C_0,\ldots,C_{p-1}$ by a unary predicate $P_0$, then choose $p$
copies $C_{i,0},\ldots,C_{i,p-1}$ between $C_i$ and $C_{(i+1)({\rm
mod}p)}$, for each $i\leq p-1$, and mark these copies by a unary
predicate $P_1$. Continuing the process with countably many
disjoint unary predicates $P_n$, $n\in\omega$, we obtain
increasing finite cardinalities for orbits of the first elements
of the marked copies of $\omega$ witnessing the equalities
(\ref{deg_inf}).

We observe the similar effect for any $n$-spherically ordered
structure $\mathcal{M}$. In particular, ${\rm acl}(A)={\rm
dcl}(A)$ for any  set $A\subseteq M$ of a cardinality $\geq n-2$,
since in such a case a $n$-spherical order is reduced to a linear
one \cite{Sud_spher_min}.}
\end{example}

{\bf Definition.} If for a theory $T$, ${\rm dcl}(A)={\rm acl}(A)$
for any set $A$ with $|A|\geq n$ then minimal such $n$ is called
the {\em ${\rm acl}$-${\rm dcl}$-difference} and denoted by ${\rm
acl}$-${\rm dcl}_{\rm dif}(T)$. If such natural $n$ does not
exists, i.e., for any $n\in\omega$ there exists a set $A$ with
$|A|\geq n$ and ${\rm acl}(A)\supset{\rm dcl}(A)$ then we put
${\rm acl}$-${\rm dcl}_{\rm dif}(T)=\infty$.

\begin{remark} \label{note_ordered}
{\rm By the definition and Example \ref{ex4} we observe that ${\rm
acl}$-${\rm dcl}_{\rm dif}(T_0)=0$ for any linearly ordered theory
$T_0$, ${\rm acl}$-${\rm dcl}_{\rm dif}(T)\leq 1$ for any
circularly ordered theory $T$, with ${\rm acl}$-${\rm dcl}_{\rm
dif}(T_1)=1$ if $T_1$ is circularly ordered with unique $1$-type
and a model with at least two elements. More generally, ${\rm
acl}$-${\rm dcl}_{\rm dif}(T)\leq n-2$ for any $n$-spherically
ordered theory $T$, where $n\geq 3$, with ${\rm acl}$-${\rm
dcl}_{\rm dif}(T_{n-2})=n-2$ if, for instance, $T_{n-2}$ is dense
$n$-spherically ordered \cite{KS_22}.

Taking a disjoint union \cite{CCMCT, Wo} of models
$\mathcal{M}_{n-2}$ of $n$-spherically ordered theories $T_{n-2}$,
with unboundedly many $n$, we obtain a structure $\mathcal{M}$
with ${\rm acl}$-${\rm dcl}_{\rm dif}({\rm
Th}(\mathcal{M}))=\infty$.}
\end{remark}

In view of Remark \ref{note_ordered} we have the following:

\begin{theorem}\label{th_ad_dif}
For any $\lambda\in\omega\cup\{\infty\}$ there is a theory
$T_\lambda$ with ${\rm acl}$-${\rm dcl}_{\rm
dif}(T_\lambda)=\lambda$.
\end{theorem}

\section{Algebraic sets and their degrees}

{\bf Definition.} Let $\mathcal{M}$ be a $L$-structure,
$A,B\subseteq M$. The set $B$ is called {\em algebraic} over $A$,
or {\em $A$-algebraic} if it is the finite set of solutions in
$\mathcal{M}$ of a $L(A)$-formula $\varphi(x,\overline{a})$, where
$\varphi(x,\overline{y})$ is a $L$-formula and $\overline{a}\in
A$, i.e., $B$ consists of algebraic elements over $A$ witnessed by
a fixed formula $\varphi(x,\overline{a})$. If $A$ is empty then an
$A$-algebraic set is called {\em algebraic}.

The set $B$ is called {\em $(A,u)$-algebraic} if it is a union of
$A$-algebraic sets.

If $B$ is a union of $A$-algebraic subsets of cardinalities at
most $n$, for $n\in\omega$, then $B$ is called {\em
$(A,u,n)$-algebraic}. If $B$ is both $A$-algebraic and
$(A,u,n)$-algebraic then $B$ is called {\em $(A,n)$-algebraic}.

The least $n$ for the $(A,u,n)$-algebraicity of $B$ is called the
{\em degree} of the $(A,u,n)$-algebraicity of $B$ and it is
denoted by ${\rm deg}_{{\rm alg},u}(B/A)$. If $B$ is
$(A,n)$-algebraic then the value ${\rm deg}_{{\rm alg},u}(B/A)$ is
defined, too, with ${\rm deg}_{{\rm alg},u}(B/A)\leq n$.

At the same time an $(A,u)$-algebraic set $B$ admits the value
${\rm deg}_{{\rm alg},u}(B/A)=\infty$, if $B$ is not represented
as a union of $A$-algebraic sets in bounded finite cardinalities.

The value ${\rm sup}\{{\rm deg}_{{\rm alg},u}(B/A)\mid B$ is
$(A,u)$-algebraic$\}$ is called the degree of the
$(A,u,n)$-algebraicity and denoted by ${\rm DEG}_{{\rm
alg},u}(A)$.

We denote by $\mathcal{A}_{A,\mathcal{M}}$ the set of all
$A$-algebraic sets in $\mathcal{M}$, by
$\mathcal{A}^u_{A,\mathcal{M}}$ the set of all $(A,u)$-algebraic
sets in $\mathcal{M}$, by $\mathcal{A}^n_{A,\mathcal{M}}$ the set
of all $(A,n)$-algebraic sets in $\mathcal{M}$, and by
$\mathcal{A}^{u,n}_{A,\mathcal{M}}$ the set of all
$(A,u,n)$-algebraic sets in $\mathcal{M}$.

We omit $A$ above if it is empty.

\medskip
Notice that the sets $\mathcal{A}_{A,\mathcal{M}}$,
$\mathcal{A}^u_{A,\mathcal{M}}$, $\mathcal{A}^n_{A,\mathcal{M}}$,
$\mathcal{A}^{u,n}_{A,\mathcal{M}}$ are preserved under elementary
extensions of $\mathcal{M}$. Therefore we omit the indexes
$\mathcal{M}$ above, denoting the sets above by $\mathcal{A}_{A}$,
$\mathcal{A}^u_{A}$, $\mathcal{A}^n_{A}$, $\mathcal{A}^{u,n}_{A}$.

\medskip
The following assertion collects some properties of algebraic
sets.

\begin{proposition}\label{alg_properties}
$1.$ For any set $A$ and $m\leq n<\omega$,
$\mathcal{A}^m_{A}\subseteq
\mathcal{A}^n_{A}\cap\mathcal{A}^{u,m}_{A}\subseteq\mathcal{A}_{A}\subseteq\mathcal{A}^u_{A}$,
$\mathcal{A}^{u,m}_{A}\subseteq
\mathcal{A}^n_{A}\cap\mathcal{A}^{u,n}_{A}\subseteq\mathcal{A}^u_{A}$.

$2.$ If ${\rm deg}_{{\rm alg},u}(B/A)$ is defined, i.e.
$B\in\mathcal{A}^u_{A}$, then $B\subseteq {\rm acl}(A)$. The
converse holds if $B$ consists of some finite orbits {\rm (}in a
saturated extension of $\mathcal{M}${\rm )} of elements under
$A$-automorphisms.

$3.$ {\rm (Monotony)} For any $B_1,B_2\in\mathcal{A}^u_{A}$ with
$B_1\subseteq B_2$, ${\rm deg}_{{\rm alg},u}(B_1/A)\leq{\rm
deg}_{{\rm alg},u}(B_2/A)$.

$4.$ {\rm (Monotony)} If $B\in\mathcal{A}^u_{A_1}$ and
$A_1\subseteq A_2$ then ${\rm deg}_{{\rm alg},u}(B/A_1)\geq{\rm
deg}_{{\rm alg},u}(B/A_2)$.

$5.$ For any set $A$, ${\rm DEG}_{{\rm alg},u}(A)={\rm deg}_{{\rm
alg},u}({\rm acl}(A)/A)$.

$6.$ For any sets $A$ and $B$, ${\rm deg}_{\rm alg}(B/A)=0$ iff
$B=\emptyset$.

$7.$ For any sets $A$ and $B$, ${\rm deg}_{\rm alg}(B/A)=1$ iff
$\emptyset\ne B\subseteq{\rm dcl}(A)$.

$8.$ If $A$ is algebraically closed then ${\rm DEG}_{{\rm
alg},u}(A)=0$ if $A$ is empty, and ${\rm DEG}_{{\rm alg},u}(A)=1$
if $A$ is nonempty.

\end{proposition}

Proof. Item 1 holds since any $m$-algebraic element over $A$ is
algebraic over $A$, and, moreover, $n$-algebraic over $A$, for any
$m\leq n$.

Item 2 is satisfied as elements of $\mathcal{A}^u_{A}$ consist of
some finite orbits of elements under $A$-automorphisms.

The monotonies are true by the definition.

Item 5 is implied by Item 2 and Monotony (Item 3).

Item 6 is obvious.

Item 7 holds by Item 6 and the property both for ${\rm deg}_{\rm
alg}(B/A)=1$ and for $B\subseteq{\rm dcl}(A)$ that $B$ is composed
as a union of $A$-definable singletons.

Item 8 is implies by $A={\rm dcl}(A)$ and Items 6 and 7.

\begin{corollary}\label{alg_alg}
If $\mathcal{M}$ is an algebra, then for any $A\subset M$ the
universe $M(A)$ of the subalgebra $\mathcal{M}(A)$ of
$\mathcal{M}$ generated by $A$ is $(A,u)$-algebraic with ${\rm
deg}_{{\rm alg},u}(M(A)/A)=1$.
\end{corollary}

Proof. Any element $a$ of $M(A)$ is represented by a term
$t(a_1,\ldots,a_n)$, where $a_1,\ldots,a_n\in A$. Therefore
$M(A)\subseteq{\rm dcl}(A)$. Thus $M(A)$ is represented as a union
of singletons $\{a\}\in\mathcal{A}^1_{A}\subseteq
\mathcal{A}^u_{A}$ whence $M(A)\in \mathcal{A}^u_{A}$ and ${\rm
deg}_{{\rm alg},u}(M(A)/A)=1$ in view of Item 7 in Proposition
\ref{alg_properties}.

\begin{remark} \label{deg_alg_val} {\rm For any set $A$, ${\rm DEG}_{{\rm
alg},u}(A)\in\omega\cup\{\infty\}$, since any $(A,u)$-algebraic
set is either represented as a union of $(A,u,n)$-algebraic sets
for some $n\in\omega$, or such $n$ does not exists for some
$(A,u)$-algebraic set $B$. In the first case we have ${\rm
DEG}_{{\rm alg},u}(A)\in\omega$, and in the second one ${\rm
DEG}_{{\rm alg},u}(A)=\infty$.

For an illustration we show that all values ${\rm DEG}_{{\rm
alg},u}(A)\in\omega\cup\{\infty\}$ are realized for appropriate
sets $A$.

Indeed, taking an equivalent relation $E$ with infinitely many
equivalence classes of cardinality $n+1$, we obtain ${\rm
DEG}_{{\rm alg},u}(\{a\})=n$ for any singleton $A=\{a\}$, where
there are exactly four $A$-algebraic sets: $\emptyset$, $\{a\}$,
$E(a)\setminus\{a\}$, $E(a)$.

Now let $E_n$, $n\in\omega$, be an ascending chain of equivalence
relations such that $E_0$-classes are singletons and each
$E_{n+1}$-class consists of exactly two $E_n$-classes. Taking an
arbitrary singleton $A=\{a\}$ we obtain ${\rm DEG}_{{\rm
alg},u}(\{a\})=\infty$, witnessed by the set
$B=\bigcup\limits_{n\in\omega}E_n(a)$, since all $E_n$-classes are
finite and have unboundedly increasing cardinalities.

We obtain a similar effect taking disjoint unary predicates $P_n$
with $|P_n|=n$, $n\in\omega$. Here ${\rm DEG}_{{\rm
alg},u}(\emptyset)=\infty$, and ${\rm DEG}_{{\rm
alg},u}(\emptyset)\in\omega$, if the language is restricted to
finitely many predicate symbols.}
\end{remark}

\section{Further variations of algebraic and definable closures and their connections}

{\bf Definition.} 1. Let $\Delta$ be a set of formulae of a theory
$T$. For a models $\mathcal{M}\models T$ and a set $A\subseteq M$
the union of sets of solutions of formulae
$\varphi(x,\overline{a})$, where
$\varphi(x,\overline{y})\in\Delta$ and $\overline{a}\in A$, such
that $\models\exists^{=n}x\:\varphi(x,\overline{a})$ for some
$n\in\omega$ (respectively,
$\models$~$\exists^{=1}x\:\varphi(x,\overline{a})$) is said to be
a \emph{$\Delta$-algebraic {\rm (}$\Delta$-definable {\rm or}
$\Delta$-definitional{\rm )} closure} of $A$. The
$\Delta$-algebraic closure of $A$ is denoted by ${\rm
acl}^{\Delta}(A)$ and its $\Delta$-definable
($\Delta$-definitional) closure, by ${\rm dcl}^\Delta(A)$.

In such a case we say that the formulae $\varphi(x,\overline{a})$
{\em witness} that $\Delta$-algebraic / $\Delta$-definable
($\Delta$-definitional) closure, and these formulae are called
{\em $\Delta$-algebraic / $\Delta$-defining}.

Any element $b\in{\rm acl}^\Delta(A)$ (respectively, $b\in{\rm
dcl}^\Delta(A)$) is called {\em $\Delta$-algebraic} ({\em
$\Delta$-definable} or {\em $\Delta$-definitional}) over $A$. If
the set $A$ is fixed or empty, we just say that $b$ is {\em
$\Delta$-algebraic}, {\em $\Delta$-definable}, or
$\Delta$-definitional.

2. If ${\rm dcl}^\Delta(A)={\rm acl}^\Delta(A)$, ${\rm
cl}_1^\Delta(A)$ denotes their common value. In such a case we say
that $A$ is {\em $\Delta$-quasi-Urbanik}.

3. If $A={\rm acl}^\Delta(A)$ (respectively, $A={\rm
dcl}^\Delta(A)$) then $A$ is called {\em $\Delta$-algebraically}
({\em $\Delta$-definably}) closed.

\medskip
We combine $\Delta$-algebraic closures and $n$-algebraic closure
as follows.

\medskip
{\bf Definition.} 1. For $n\in\omega\setminus\{0\}$, a set
$A\subseteq M$ and a set $\Delta$ of formulae an element $b$ is
called {\em $(\Delta,n)$-algebraic} over $A$, if $b\in{\rm
acl}^\Delta(A)$ and it is witnessed by a formula
$\varphi(x,\overline{a})$, for some $\overline{a}\in A$ and
$\varphi(x,\overline{y})\in\Delta$, with at most $n$ solutions.

2. The set of all $(\Delta,n)$-algebraic elements over $A$ is
denoted by ${\rm acl}^\Delta_n(A)$.

3. If $A={\rm acl}^\Delta_n(A)$ then $A$ is called {\em
$(\Delta,n)$-algebraically} closed.

4. If ${\rm acl}^\Delta(A) = {\rm acl}^\Delta_n(A)$ then minimal
such $n$ is called the {\em $\Delta$-degree of algebraization}
over the set $A$ and it is denoted by ${\rm deg}^\Delta_{\rm
acl}(A)$. If that $n$ does not exist then we put ${\rm
deg}^\Delta_{\rm acl}(A)=\infty$. The supremum of values ${\rm
deg}^\Delta_{\rm acl}(A)$ with respect to all sets $A$ of given
theory $T$ is denoted by ${\rm deg}^\Delta_{\rm acl}(T)$ and
called the {\em $\Delta$-degree of algebraization} of the theory
$T$.

5. A theory $T$ with ${\rm deg}^\Delta_{\rm acl}(T)=1$, i.e., with
quasi-Urbanik sets $A$ of $T$ only, is called {\em
$\Delta$-quasi-Urbanik}, and the models $\mathcal{M}$ of $T$ are
{\em $\Delta$-quasi-Urbanik}, too.

\medskip
The following remark collects some obvious properties related to
the operators ${\rm acl}^\Delta$ and ${\rm acl}^\Delta_n$.

\begin{remark} \label{Delta-acl0} {\rm
$1.$ {\rm (Monotony)} If $m\leq n$, $\Delta_1\subseteq\Delta_2$,
and $A_1\subseteq A_2$ then ${\rm
acl}^{\Delta_1}_m(A_1)\subseteq{\rm acl}^{\Delta_2}_n(A_2)$.

$2.$ {\rm ($\Delta$-reduction)} For any set $A$, ${\rm
acl}(A)={\rm acl}^\Delta(A)$ {\rm (}respectively, ${\rm
acl}_n(A)={\rm acl}_n^\Delta(A)$, where
$n\in\omega\setminus\{0\}${\rm )} iff for any $a\in{\rm acl}(A)$
{\rm (}$a\in{\rm acl}_n(A)${\rm )} it is witnessed by a formula in
$\Delta$.}
\end{remark}

\medskip
{\bf Definition} \cite{ZPS}). A theory $T$ is said to be
\emph{$\Delta$-based}\index{Theory!$\Delta$-based}, where $\Delta$
is some set of~formulas without parameters, if any formula of $T$
is equivalent in $T$ to a~Boolean combination of formulae
of~$\Delta$.

For $\Delta$-based theories $T$, it is also said that $T$ has {\em
quantifier elimination}\index{Elimination of quantifiers} or {\em
quantifier reduction}\index{Reduction of quantifiers} up to
$\Delta$.

\medskip
Let $\Delta$ be a set of formulae of a theory $T$.

\begin{remark} \label{Delta-based} {\rm
1. For the theory $T$, the operators ${\rm acl}^\Delta$ and ${\rm
acl}^\Delta_n$ coincide iff finite definable sets $A$ for formulae
$\varphi(x,\overline{y})\in\Delta$ with
$A=\varphi(\mathcal{M},\overline{a})$, $\mathcal{M}\models T$ are
covered by definable sets $A_i$ for formulae
$\varphi_i(x,\overline{y})\in\Delta$ with
$A_i=\varphi_i(\mathcal{M},\overline{a})$ having at most $n$
elements.

2. If a theory $T$ is $\Delta$-based then both ${\rm
acl}^{\Delta'}={\rm acl}$ and ${\rm acl}^{\Delta'}_n={\rm acl}_n$,
where
$$\Delta'=\left\{\bigwedge\limits_i\varphi_i^{\delta_i}\mid\varphi_i\in\Delta,\delta_i\in\{0,1\}\right\}.$$
}
\end{remark}

\begin{proposition}\label{Delta-acl} $1.$ The inclusion $A\subseteq{\rm
acl}^\Delta(A)$ {\rm ({\em respectively,} $A\subseteq{\rm
acl}^\Delta_n(A)$)} holds for any set $A$ if and only if $\Delta$
is {\rm ${\rm acl}$-reflexive} {\rm ($({\rm acl},n)$-reflexive)},
i.e., for any element $a$ there is $\varphi(x,y)\in\Delta$ with
$\models\varphi(a,a)\wedge\exists^{\leq m}x\,\varphi(x,a)$ for
some $m\in\omega$ {\rm ({\em for some $m\leq n$, where
$n\in\omega\setminus\{0\}$})}.

$2.$ The inclusion ${\rm acl}^\Delta({\rm
acl}^\Delta(A))\subseteq{\rm acl}^\Delta(A)$ {\rm ({\em
respectively,} ${\rm acl}^\Delta_n({\rm
acl}^\Delta_n(A))\subseteq{\rm acl}^\Delta_n(A)$)} holds for any
set $A$ if and only if $\Delta$ is {\rm ${\rm acl}$-transitive}
{\rm ($({\rm acl},n)$-transitive)}, i.e., for any tuples
$\overline{a}\in A$, $\overline{b}\in{\rm acl}^\Delta(A)$ {\rm
($\overline{b}\in{\rm acl}^\Delta_n(A)$)} and a formula
$\varphi(x,\overline{z})\in\Delta$ with
$\models\varphi(c,\overline{b})\wedge\exists^{<\omega}x\varphi(x,\overline{b})$
$(\models\varphi(c,\overline{b})\wedge\exists^{\leq
n}x\varphi(x,\overline{b}))$ there is a formula
$\psi(x,\overline{y})\in\Delta$ with
$\models\psi(c,\overline{a})\wedge\exists^{<\omega}x\psi(x,\overline{a})$
$(\models\psi(c,\overline{a}\wedge\exists^{\leq
n}x\psi(x,\overline{a})$, where $n\in\omega\setminus\{0\})$.

$3.$ The operators ${\rm acl}^\Delta$ and ${\rm acl}^\Delta_n$
satisfy Finite character.

$4.$ The operator ${\rm acl}$ {\rm (}${\rm acl}_n${\rm )} has
Exchange property iff any operator ${\rm acl}_n^\Delta$ {\rm
(}${\rm acl}_n^\Delta${\rm )} is extensible till some ${\rm
acl}^{\Delta'}$ {\rm (}${\rm acl}^{\Delta'}_n${\rm )} with
$\Delta\subseteq\Delta'$ and Exchange property. Here either ${\rm
acl}^{\Delta'}$ {\rm (}${\rm acl}^{\Delta'}_n$ with given $n${\rm
)} has a proper extension or does not have proper extensions
depending on the given theory.
\end{proposition}

Proof is obvious using the definition of operators ${\rm
acl}^\Delta$ and ${\rm acl}^\Delta_n$.

\medskip
Proposition \ref{Delta-acl} characterizes basic properties of
closure operators. If the operator ${\rm acl}^\Delta$
(respectively, ${\rm acl}^\Delta_n$) satisfies the conditions
described in Items 1 and 2 of Proposition \ref{Delta-acl} then it
is called {\em regular}. If, additionally, the condition in Item 4
holds, then that operator is called {\em pregeometric}.

Proposition \ref{Delta-acl} immediately implies:

\begin{corollary}\label{cor_pregeom}
The operator ${\rm acl}^\Delta$ {\rm ({\em respectively,} ${\rm
acl}^\Delta_n$)} defines a pregeometry on the universe of given
structure iff it is pregeometric, i.e., it is regular and
satisfies Exchange property.
\end{corollary}

\begin{remark} \label{tran_clos}{\rm
Any operator ${\rm acl}^\Delta$ is extensible till transitive one
${\rm acl}^{\Delta'}$ by the transitive closure forming $\Delta'$
by adding to $\Delta$ all formulae
\begin{equation}\label{eq_acl2}\exists
x_1,\ldots,x_n\left(\psi(x,x_1,\ldots,x_n)\wedge\bigwedge\limits_{i=1}^n\varphi_i(x_i,\overline{y})\right)\end{equation}
as in (\ref{eq_acl1}), where $\varphi_i$ and $\psi$ either belong
to $\Delta$ or were added before, i.e., $\Delta'$ is the set of
formulae obtained from $\Delta$ in the calculus with the rules:
$$
\displaystyle\frac{\varphi_1,\ldots,\varphi_n,\psi}{\exists
x_1,\ldots,x_n\left(\psi(x,x_1,\ldots,x_n)\wedge\bigwedge\limits_{i=1}^n\varphi_i(x_i,\overline{y})\right)},\mbox{\
\ for } n\in\omega.
$$

Using Monotony in Remark \ref{Delta-acl0} we have ${\rm
acl}^\Delta(A)\subseteq{\rm acl}^{\Delta'}(A)$ for any set $A$ in
the given structure. Taking $\Delta''=\Delta'\cup\{x\approx y\}$
we obtain a regular operator ${\rm acl}^{\Delta''}$.

We have a similar regular operator ${\rm acl}^{\Delta''}_1$
starting with the operator ${\rm acl}^{\Delta}_1$.

In contrast with ${\rm acl}^\Delta$ and  ${\rm acl}^{\Delta}_1$
the transitivity of ${\rm acl}^{\Delta''}_n$, for $n\geq 2$, can
fail under the transitive closure above, since compositions of
$n$-algebraic formulae may not be $n$-algebraic.}
\end{remark}

In view of properties above a hierarchy of regular and
pregeometric algebraic closure operators ${\rm acl}^\Delta$, ${\rm
acl}^\Delta_n$ arises depending on chosen sets $\Delta$ of
formulae and bounds $n$ of cardinalities for sets of solutions of
algebraic formulae in $\Delta$. This hierarchy includes various
degrees ${\rm acl}$-${\rm dcl}_{\rm dif}(T)$ of algebraization.

\medskip
{\bf Definition.} Two sets $\Delta_1$ and $\Delta_2$ of formulae
are called {\em ${\rm acl}$-equivalent} (respectively, ${\rm
acl}_n$-equivalent), denoted by $\Delta_1\sim\Delta_2$
($\Delta_1\sim_n\Delta_2$), if ${\rm acl}^{\Delta_1}={\rm
acl}^{\Delta_2}$ (${\rm acl}^{\Delta_1}_n={\rm
acl}^{\Delta_2}_n$).

\begin{remark} \label{inter}
{\rm Equivalent sets preserve the regularity and pregeometricity.
Moreover, equivalence classes $\sim\!\!(\Delta_1)$ and
$\sim\!\!(\Delta_2)$ (respectively, $\sim_n\!\!(\Delta_1)$ and
$\sim_n\!\!(\Delta_2)$) allow to choose {\em coordinated}
representatives $\Delta'_1\in\,\,\sim\!\!(\Delta_1)$ and
$\Delta'_2\in\,\,\sim\!\!(\Delta_2)$ (respectively,
$\Delta'_1\in\,\,\sim_n\!\!(\Delta_1)$ and
$\Delta'_2\in\,\,\sim_n\!\!(\Delta_2)$) such that for regular
${\rm acl}^{\Delta'_1}$ and ${\rm acl}^{\Delta'_2}$, ${\rm
acl}^{\Delta'_1}\cap\,\,{\rm acl}^{\Delta'_2}$ (${\rm
acl}_n^{\Delta'_1}\cap\,\,{\rm acl}_n^{\Delta'_2}$) is regular,
too, and ${\rm acl}^{\Delta'_1}\cap\,\,{\rm acl}^{\Delta'_2}={\rm
acl}^{\Delta'_1\cap\Delta'_2}$ (${\rm acl}_n^{\Delta'_1}\cap{\rm
acl}_n^{\Delta'_2}={\rm acl}_n^{\Delta'_1\cap\Delta'_2}$). Indeed,
it suffices to add to $\Delta'_1$ and to $\Delta'_2$ same formulae
witnessing common finite definable sets with respect to ${\rm
acl}^{\Delta_1}$ and ${\rm acl}^{\Delta_2}$ (${\rm
acl}_n^{\Delta_1}$ and ${\rm acl}_n^{\Delta_2}$).

Thus, the binary operations $\wedge$ and $\wedge_n$ of
intersection arise mapping the pairs $({\rm acl}^{\Delta_1}, {\rm
acl}^{\Delta_2})$ (respectively $({\rm acl}_n^{\Delta_1},{\rm
acl}_n^{\Delta_2})$) to ${\rm acl}^{\Delta'_1\cap\Delta'_2}$
(${\rm acl}_n^{\Delta'_1\cap\Delta'_2}$). }
\end{remark}

For a structure $\mathcal{M}$ we define derived structures
$\mathcal{SL}_{\rm acl}(\mathcal{M})$ and $\mathcal{SL}_{{\rm
acl}_n}(\mathcal{M})$, $n\in\omega$, in the following way. The
universe ${\rm SL}_{\rm acl}(\mathcal{M})$ (respectively, ${\rm
SL}_{{\rm acl}_n}(\mathcal{M})$) of $\mathcal{SL}_{\rm
acl}(\mathcal{M})$ ($\mathcal{SL}_{{\rm acl}_n}(\mathcal{M})$)
consists of all regular operators ${\rm acl}^\Delta$ and ${\rm
acl}^{\Delta}_n$ for $\mathcal{M}$, and the language consists of
the symbol $\wedge$ ($\wedge_n$) for the operation of
intersection.

If  $\mathcal{M}$ is sufficiently saturated then the operators
${\rm acl}^{\Delta}$ and ${\rm acl}^{\Delta}_n$ admit syntactical
descriptions by families of formulae, and it induces structures
$\mathcal{SL}_{\rm acl}(T)$ and $\mathcal{SL}_{{\rm acl}_n}(T)$,
for the theory $T={\rm Th}(\mathcal{M})$, which are isomorphic to
$\mathcal{SL}_{\rm acl}(\mathcal{M})$ and $\mathcal{SL}_{{\rm
acl}_n}(\mathcal{M})$, respectively.

\begin{proposition}\label{semilat}
$1.$ Any structure $\mathcal{SL}_{\rm acl}(T)$ is a lower
semilattice with the least and the greatest elements.

$2.$ Any structure $\mathcal{SL}_{{\rm acl}_n}(T)$ is a lower
semilattice with the least element. It has the greatest element if
${\rm acl}_n$ is regular.
\end{proposition}

Proof. By the definition both intersections ${\rm
acl}^{\Delta'_1}\cap\,\,{\rm acl}^{\Delta'_2}={\rm
acl}^{\Delta'_1\cap\Delta'_2}$ and ${\rm
acl}_n^{\Delta'_1}\cap{\rm acl}_n^{\Delta'_2}={\rm
acl}_n^{\Delta'_1\cap\Delta'_2}$ in Remark \ref{inter} are infimum
for $({\rm acl}^{\Delta_1}$, ${\rm acl}^{\Delta_2})$ and $({\rm
acl}_n^{\Delta_1}$, ${\rm acl}_n^{\Delta_2})$, respectively. The
operator ${\rm acl}^{\{x\approx y\}}$ is the least element both
for $\mathcal{SL}_{\rm acl}(T)$ and $\mathcal{SL}_{{\rm
acl}_n}(T)$. The operator ${\rm acl}$ is the greatest element for
$\mathcal{SL}_{\rm acl}(T)$. And ${\rm acl}_n$ is the greatest
element for $\mathcal{SL}_{\rm acl_n}(T)$ if ${\rm acl}_n$ is
regular.

\medskip
Since ${\rm acl}_1$ is always regular, Proposition \ref{semilat}
implies the following:

\begin{corollary}
Any structure $\mathcal{SL}_{{\rm acl}_1}(T)$ is a lower
semilattice with the least and the greatest elements.
\end{corollary}

The following example illustrates that $\mathcal{SL}_{{\rm
acl}_n}(T)$ may not have the greatest element.

\begin{example}\label{ex_ng}
{\rm Let $\mathcal{H}=\langle M,Z,E_1,E_2\rangle$,
$Z\subseteq\mathcal{P}(M)$, be a hypergraph \cite{EmMe} with
colored hyperedges, being equivalence classes of equivalence
relations $E_1$ and $E_2$, and satisfying the following
conditions:

i) each hyperedge consists of three elements;

ii) any element belongs to exactly two hyperedges, and one of
these hyperedges is marked by $E_1$ and another one by $E_2$;

iii) the hypergraph does not have cycles.

For the structure $\mathcal{M}=\langle M,E_1,E_2\rangle$ we have
regular $\subseteq$-incomparable operators ${\rm acl}^{\{x\approx
y,E_i(x,y)\}}_2$, $i=1,2$, which do not have proper regular
extensions, since the operator ${\rm acl}_2$ on $M$ is not
transitive. Indeed, taking an arbitrary element $a\in M$ we have a
$2$-algebraic formulae $E_1(x,a)\wedge\neg (x\approx a)$ and
$E_1(x,a)\wedge\neg (x\approx a)$. But their composition is not
$2$-algebraic having $4$ solutions.

Thus the semilattice $\mathcal{SL}_{{\rm acl}_2}(\mathcal{M})$
consists of three operators: the least one ${\rm acl}^{\{x\approx
y\}}_2$ and two its $\subseteq$-incomparable extensions ${\rm
acl}^{\{x\approx y,E_i(x,y)\}}_2$, $i=1,2$. Moreover, since
cardinalities of isolating algebraic formulae $\varphi(x,a)$ are
unbounded, for iterated compositions of $E_1(x,y)$ and $E_2(x,y)$,
operators ${\rm acl}_n$, for $n\geq 3$, are not regular, too,
implying that the semilattices $\mathcal{SL}_{{\rm
acl}_n}(\mathcal{M})$ also consist of three operators and do not
have the greatest elements. Adding to $\mathcal{SL}_{{\rm
acl}_n}(\mathcal{M})$ the operator ${\rm acl}$ we obtain the
$4$-element semilattice $\mathcal{SL}_{{\rm acl}}(\mathcal{M})$
forming a Boolean algebra. Notice also that $\mathcal{SL}_{{\rm
acl}_1}(\mathcal{M})$ is a singleton consisting of the operator
${\rm acl}^{\{x\approx y\}}_1$.

We obtain the similar effect increasing the cardinality of
$E_i$-classes till natural $m\geq 4$. At the same time, if $m=2$,
then $\mathcal{SL}_{{\rm acl}_1}(\mathcal{M})$ forms a $4$-element
lattice with the least and the greatest elements.

Thus for the theories $T_m$ of structures $\mathcal{M}_m=\langle
M_m,E_1,E_2\rangle$ consisting of $m$-element hyperedges marked by
$E_1$ and $E_2$ as above, either ${\rm deg}_{\rm acl}(T_m)=1$, if
$m=1$ or $m=2$, or ${\rm deg}_{\rm acl}(T_m)=\infty$, if $m\geq
3$.}
\end{example}

Now we define expansions of the structures $\mathcal{SL}_{\rm
acl}(\mathcal{M})$ and $\mathcal{SL}_{{\rm acl}_n}(\mathcal{M})$,
$n\in\omega$, by the operation of union $\vee$. For the regular
operators ${\rm acl}^{\Delta_1}$ and ${\rm acl}^{\Delta_2}$
(respectively, ${\rm acl}_n^{\Delta_1}$ and ${\rm
acl}_n^{\Delta_2}$) we put ${\rm acl}^{\Delta_1}\vee{\rm
acl}^{\Delta_2}={\rm acl}^{\Delta}$ (${\rm
acl}_n^{\Delta_1}\vee{\rm acl}_n^{\Delta_2}={\rm
acl}_n^{\Delta}$), where $\Delta$ is the transitive closure of
$\Delta_1\cup\Delta_2$. We denote these expansions of
$\mathcal{SL}_{\rm acl}(\mathcal{M})$ and $\mathcal{SL}_{{\rm
acl}_n}(\mathcal{M})$ by $\mathcal{L}_{\rm acl}(\mathcal{M})$ and
$\mathcal{L}_{{\rm acl}_n}(\mathcal{M})$, respectively. Again
choosing $\mathcal{M}$ sufficiently saturated we obtain structures
$\mathcal{L}_{\rm acl}(T)$ and $\mathcal{L}_{{\rm acl}_n}(T)$ for
$T={\rm Th}(\mathcal{M})$.

Since the union $\vee$ and the intersection $\wedge$ correspond to
appropriate union and intersection of $\Delta_i$, in view of
Proposition \ref{semilat} we have the following:

\begin{theorem}\label{lat}
The structure $\mathcal{L}_{\rm acl}(T)$ is a distributive lattice
with the least and the greatest elements.
\end{theorem}

Notice that the structures $\mathcal{L}_{{\rm
acl}_n}(\mathcal{M})$ can be a distributive lattices, too, taking,
for instance, finite $\mathcal{M}$ and $n\geq |M|$. At the same
time, as Example \ref{ex_ng} shows, in general case the structures
$\mathcal{L}_{{\rm acl}_n}(\mathcal{M})$ may be not lattices.

\medskip
Recall that for a lattice $\mathcal{L}$, the lattice {\em height}
({\em width}) is the supremum of cardinalities for (anti)chains.
We denote these characteristics by $h(\mathcal{L})$ and
$w(\mathcal{L})$, respectively,

The following theorem shows that height and width of
$\mathcal{L}_{\rm acl}(T)$ can be unbounded.

\begin{theorem}\label{lat_hw}
For any cardinality $\lambda>0$ there is a lattice
$\mathcal{L}_{\rm acl}(T)$ with $h(\mathcal{L}_{\rm
acl}(T))=\lambda$ and $w(\mathcal{L}_{\rm acl}(T))=\lambda$.
\end{theorem}

Proof. The case $\lambda=1$ is realized by an arbitrary theory $T$
of a singleton. The case $\lambda=2$ is realized in Example
\ref{ex_ng}. Now for $\lambda\geq 3$ we extend a structure
$\mathcal{N}\equiv\mathcal{M}=\langle M,E_1,E_2\rangle$ for
Example \ref{ex_ng} till a structure $\mathcal{N}'$ with $\lambda$
equivalence relations $E_i$, $i<\lambda$, satisfying the following
conditions:

i) each $E_i$-class contains 3 elements;

ii) for every $E_i$-class $X$ and $E_j$-class $Y$, where $i\ne j$,
$|X\cap Y|\leq 1$;

iii) the hypergraph $\langle M,Z\rangle$, where $Z$ is the set of
all $E_i$-classes, $i<\lambda$.

For the theory $T={\rm Th}(\mathcal{N}')$ the lattice
$\mathcal{L}_{\rm acl}(T)$ has $\lambda$ atoms  ${\rm
acl}^{\{x\approx y,E_i(x,y)\}}$, $i<\lambda$, witnessing
$w(\mathcal{L}_{\rm acl}(T))=\lambda$. Collecting formulae
$E_i(x,y)$ we obtain a chain of $\lambda$ operators ${\rm
acl}^{\Delta_\mu}$, where $\Delta_\mu={\{x\approx y\}\cup\{E_i\mid
i<\mu\}}$, $\mu\leq\lambda$, witnessing that  $h(\mathcal{L}_{\rm
acl}(T))=\lambda$.

\section{Conclusion}

We studied variations of algebraic closure, properties and
possibilities of these variations, of degree of algebraization,
and of difference between definable and algebraic closures. These
possibilities are illustrated by a series of examples. Algebraic
sets and their degrees are studied. Hierarchy of operators of
algebraic closures relative sets is considered, semilattices and
lattices for families of these operators are introduced and some
characteristics of these structures are described.

It would be natural to describe degrees of algebraization and
difference between definable and algebraic closures for natural
classes of structures and their theories. It would be also
interesting to describe possibilities of Hasse diagrams for
semilattices and lattices for families of operators of algebraic
closure.

Sergey V. Sudoplatov

Sobolev Institute of Mathematics,

Novosibirsk State Technical University,

Novosibirsk, Russia

E-mail: sudoplat@math.nsc.ru

\end{document}